\documentclass[12pt]{amsart}

\usepackage{amsmath,amsfonts,amssymb,amscd}
\begin{document}
\renewcommand{\thefootnote}{\fnsymbol{footnote}}
\pagestyle{plain}

\title{An affine sphere equation associated to Einstein toric surfaces}
\author{Toshiki Mabuchi}
\address{Department of Mathematics, Graduate School of Science, Osaka
University,  Toyonaka, Osaka, 560-0043 Japan
}
\maketitle

\abstract
\noindent As seen in the works of Calabi [1], Cheng-Yau [3] and Loftin [8], affine sphere equations 
have a close relationship with K\"ahler-Einstein metrics. 
The main purpose of this note is to show that an equation analogous to those of hyperbolic affine spheres
arises naturally from K\"ahler-Einstein metrics on Einstein toric surfaces. 
The case for the remaining toric surfaces with K\"ahler-Ricci solitons will also be discussed.
\endabstract


\section{Introduction}

In this note, we consider strictly convex bounded $C^2$-domains $\Omega_i$, $i=1,2,3$, 
in $\Bbb R^2:=\{(s,t)\}$ defined as the intersections
$$
\Omega_1 := \cap_{i=1}^3 \Omega_{i,1},\;\;\;
\Omega_2 := \cap_{i=1}^4 \Omega_{i,2},\;\;\; 
\Omega_3 := \cap_{i=1}^6 \Omega_{i,3},
$$
where $\Omega_{i,j}$ are the open subsets  
$\{\,(s,t)\in \Bbb R^2\;\big |\, \rho_{ij}(s,t) >0\,\}$ of $\Bbb R^2$
with the functions $\rho_{ij} = \rho_{ij}(s,t)$ defined by
$$
\begin{cases}
 & \rho_{11}:= \frac{5}{8}-\frac{3}{2}\, (s-t)^2 -
\frac{1}{2}(s+t),\\
 & \rho_{21}:= s - \frac{3}{2}\, (t + \frac{1}{6})^2  +\frac{2}{3},\\
 & \rho_{31}:= t -\frac{3}{2}\,(s + \frac{1}{6})^2+
\frac{2}{3}.
\end{cases}
$$
$$
\begin{cases}
&\rho_{12}:= -t-\frac{3}{2}\,s^2+
\frac{1}{2},\\
&\rho_{22}:= -s-\frac{3}{2}\,t^2+\frac{1}{2},\qquad\quad\;\;\;\quad\\
&\rho_{32}:= t-\frac{3}{2}\,s^2+\frac{1}{2},\\
&\rho_{42}:= s-\frac{3}{2}\,t^2+\frac{1}{2},
\end{cases}
$$
$$
\begin{cases}
& \rho_{13}:= -t-\frac{3}{2}\,(s-t+\frac{1}{6})^2+\frac{1}{3},\\ 
&\rho_{23}:= -s-\frac{3}{2}\,(t - \frac{1}{6})^2+\frac{1}{3},\\ 
&\rho_{33}:= t-\frac{3}{2}\,(s+\frac{1}{6})^2+\frac{1}{3},\\ 
&\rho_{43}:= t-\frac{3}{2}\,(-s+t+\frac{1}{6})^2+\frac{1}{3},\\ 
&\rho_{53}:=  s-\frac{3}{2}\,(t+\frac{1}{6})^2+\frac{1}{3},\\ 
&\rho_{63}:= - t-\frac{3}{2}\,(s-\frac{1}{6})^2+\frac{1}{3}.
\end{cases}
$$
Then each boundary curve $\partial\Omega_j:= \bar{\Omega}_j\setminus\Omega_j$ 
in $\Bbb R^2$ is not only 
$C^2$
but also piecewise quadratic. 
By a theorem of Cheng and Yau [2], 
there exists a unique convex negative solution $\varphi = \varphi (s,t) \in C^{\infty}(\Omega_j) \cap C^0(\bar{\Omega}_j)$  for
$$
\begin{cases}
&(-\varphi )^{2+k}\det\operatorname{Hess}(\varphi )= 1\quad\text{on $\Omega_j$},\\
	&\; \varphi_{|\partial\Omega_j} = 0,
\end{cases}\leqno{(1.1)}
$$
where the convex negativity of the solution $\varphi$ for (1.1) means that $\varphi <0$ 
on $\Omega_j$ and that the Hessian matrix
$$
\operatorname{Hess}(\varphi ):= \begin{pmatrix} 
\varphi_{ss} & \varphi_{st}\\
\varphi_{st} & \varphi_{tt}
\end{pmatrix}
$$
is positive definite everywhere on $\Omega_j$. Then for $k =2$, the equation (1.1) 
is known as the equation for hyperbolic affine spheres. 
In this note, we assume $k =1/2$,
and consider the solution $\varphi$ of (1.1).
Then by setting
$$
\psi := -\left(\frac{2}{3}\right)^{\frac{2}{3}} 
(-\varphi)^{\frac{3}{2}}\; \in\;\;  C^{\infty}(\Omega_j) 
\cap C^0(\bar{\Omega}_j),
\leqno{(1.2)}
$$
we can rewrite (1.1) in the form
$$
\begin{cases}
\;&\vmatrix
\psi_{ss} & \psi_{st} & \psi_s\\
\psi_{st} & \psi_{tt} & \psi_t\\
\psi_s    & \psi_t    & 3\psi
\endvmatrix = -3\qquad \text{on $\Omega_j$},\\
\; &\;\psi_{|\partial\Omega_j} = 0.
\end{cases} \leqno{(1.3)}
$$
In a neighborhood of $\partial\Omega_j$ in $\Bbb R^2$, we fix a $C^2$-function $\rho\,$ defining $\partial\Omega_j$ such that the 1-form $d\rho$ coincides with $d\rho_{ij}$ 
when restricted to 
$\partial\Omega_j \cap \{\rho_{ij} =0\}$ for all $i$. 
Then $\psi$ is expressible as
$- \rho - f\rho^2 + \text{higher order terms in $\rho$}$.
By abuse of terminology, we call the restriction 
$$
P (\Omega_j)\; := \; f_{|\partial\Omega_j} 
$$ 
the ``Fubini-Pick invariant'' of the domain $\Omega_j$ (cf. [11]). We now put
$X_1:= \Bbb P^2(\Bbb C)$, $X_2:= \Bbb P^1(\Bbb C)\times \Bbb P^1(\Bbb C)$,
$X_3:= \Bbb P^2(\Bbb C)\#3\bar{\Bbb P}^2(\Bbb C)$. 
Since $X_j$, $j=1,2,3$, are toric surfaces, we have natural torus 
embeddings  
$$
T :=(\Bbb C^*)^2 \hookrightarrow X_j.
$$
In this note, by setting $k:=1/2$, we shall show that the equation for
K\"ahler-Einstein metrics on $X_j$, $j =1,2,3$,
has a reduction to (1.1) above, 
where $P(\Omega_j)$ is uniquely determined by 
the pullback to $X_j \setminus T$ of the K\"ahler-Einstein form on $X_j$.
Moreover, from the data $P(\omega_j)$, we can explicitly
describe the K\"ahler-Einstein metric on $X_j$.

\section{Reduction to (1.1)}

For toric surfaces $X_j$ in the introduction,
we consider a $K$-invariant K\"ahler-Einstein form $\omega$ on $X_j$ 
in the class $2\pi c_1(X_j)$ (cf. [12],\,[13],\,[14]),
 where $K:= S^1 \times S^1$ denotes the maximal compact subgroup
 of the algebraic torus $T = \Bbb C^*\times \Bbb C^*$. 
In view of the torus embedding
$$
T =  \{ (z_1, z_2)\in \Bbb C^*\times \Bbb C^* \} \hookrightarrow X_j,
$$
we can regard $(z_1,z_2)$ as a system of holomorphic local coordinates 
on the Zariski open dense subset $T$ of $X_j$. Then the restriction to $T$ 
of the volume form $\omega^2$ on $X_j$ is written as
$$
{\omega^2}_{|T}\; =\; 2\, e^{-h} \left (\sqrt{-1}\,\frac{dz_1\wedge d\bar{z}_1}{|z_1|^2} 
\right )\wedge \left (\sqrt{-1}\,\frac{dz_2\wedge d\bar{z}_2}{|z_2|^2} 
\right )
$$
for some $K$-invariant function 
$h \in C^{\infty}(T)_{\Bbb R}$ on $T$.
Define $K$-invariant functions $x$, $y \in C^{\infty}(T)_{\Bbb R}$ on $T$ 
by
$$
e^x = |z_1|^2\quad\text{and} \quad e^y = |z_2|^2,
$$
and these are seen as real-valued independent variables
with ranges $-\infty <x < +\infty$ and $-\infty <y < +\infty$.
In particular, $h$ is regarded as a smooth function  
$$
h= h(x,y) \; \in \; C^{\infty}(\Bbb R^2)
$$
on $\Bbb R^2 = \{(x,y)\}$. 
By setting $\operatorname{Hess}(h):= 
\begin{pmatrix} 
h_{xx} & h_{xy}\\
h_{xy} & h_{yy}
\end{pmatrix}$, we see that
$$
{\operatorname{Ric}(\omega )^2}_{|T}\; =\;2\,\det \operatorname{Hess}(h) 
\left (\sqrt{-1}\,\frac{dz_1\wedge d\bar{z}_1}{|z_1|^2} 
\right )\wedge \left (\sqrt{-1}\,\frac{dz_2\wedge d\bar{z}_2}{|z_2|^2} 
\right ).
$$
Since $\omega$ is a K\"ahler-Einstein form, we have 
$\operatorname{Ric}(\omega ) = \omega$, and hence
$$
\det\operatorname{Hess}(h) = e^{-h}.
\leqno{(2.1)}
$$
Let $\square_j := \bar{\square}_j \setminus \partial {\square}_j$ be the interior of 
$\bar{\square}_j$, where
$\bar{\square}_j$ is the compact convex polygon in $\Bbb R^2 =\{(u,v)\}$ defined by
$$
\bar{\square}_ j = 
\begin{cases} 
\;\;\{ (u,v) \in \Bbb R^2\,;\, u+v \leq 1, u \geq -1, v \geq -1\},&\;\; j =1,\\
\;\;\{ (u,v) \in \Bbb R^2\,;\, |u| \leq 1, |v|\leq 1\}, &\;\; j=2,\\
\;\;\{ (u,v) \in \Bbb R^2\,;\, |u+v|\leq 1,|u|\leq 1, |v|\leq 1\}, &\;\; j=3.
\end{cases}
$$
Put 
$u:= h_x = \partial h/\partial x$ and $v:= h_y = \partial h/\partial y$.
Since the moment map
sending each $(x,y)\in \Bbb R^2$ to $(u(x,y), v(x,y))\in \square_j$
defines a diffeomorphism between $\Bbb R^2$ and $\square_j$, every function on $\Bbb R^2 = \{(x,y)\}$ is naturally regarded as a function 
on $\square_j =\{(u,v)\}$ via this moment map, and vice versa.
We now consider the Legendre transform $h^*:= xu + yv - h$. Then 
$h^*=h^*(u,v)$ 
regarded as a function in $u$ and $v$ satisfies 
$$
x = h^*_u \quad \text{and}\quad y=h^*_v,
$$
where $h^*_u := \partial h^*/\partial u$ and 
$h^*_v := \partial h^*/\partial v$. Then by (2.1) and
$$
\begin{pmatrix}
h^*_{uu} & h^*_{uv}\\
h^*_{uv} & h^*_{vv}
\end{pmatrix}\;
= \;\begin{pmatrix}
\frac{\partial x}{\partial u} & \frac{\partial x}{\partial v} \\
\frac{\partial y}{\partial u}  & \frac{\partial y}{\partial v} 
\end{pmatrix}\;
=\;\begin{pmatrix}
\frac{\partial u}{\partial x} & \frac{\partial u}{\partial y} \\
\frac{\partial v}{\partial x}  & \frac{\partial v}{\partial y} 
\end{pmatrix}^{-1}\;
=\;\begin{pmatrix}
h_{xx} & h_{xy}\\
h_{xy} & h_{yy}
\end{pmatrix}^{-1},
\leqno{(2.2)}
$$
we have $e^h h_{xx} = h^*_{vv}$, $e^h h_{yy}= h^*_{uu}$ and 
$e^h h_{xy} = -h^*_{uv}$. Now by $h^*_{vvu} = h^*_{uvv}$ and 
$h^*_{uuv} = h^*_{uvu}$, we see that
$$
\frac{\partial (e^h h_{xx})}{\partial u} 
= \frac{\partial (-e^h h_{xy})}{\partial v}
\quad\text{and}\quad
 \frac{\partial (e^h h_{yy})}{\partial v} 
= \frac{\partial (-e^h h_{xy})}{\partial u}.
\leqno{(2.3)}
$$
From the first equality of (2.3), we obtain
$$
\frac{\partial h}{\partial u}\,h_{xx} + \frac{\partial h}{\partial v}\,h_{xy}
+ \frac{\partial h_{xx}}{\partial u} + \frac{\partial h_{xy}}{\partial v} =0.
$$
Hence, together with $\frac{\partial h}{\partial u}= \frac{\partial x}{\partial u}h_x
 +  \frac{\partial y}{\partial u}h_y$ and 
$\frac{\partial h}{\partial v}= \frac{\partial x}{\partial v}h_x 
  +  \frac{\partial y}{\partial v}h_y$, it now follows that
$$
(h^*_{uu} u + h^*_{uv}v)h_{xx} + (h^*_{uv} u + h^*_{vv}v)h_{xy}+
\frac{\partial h_{xx}}{\partial u} + \frac{\partial h_{xy}}{\partial v}
=0,
$$
where we used (2.2) and the definitions of $u$ and $v$.
Again by (2.2), $h^*_{uu}h_{xx}+h^*_{uv}h_{xy} =1$ and 
$h^*_{uv}h_{xx}+h^*_{vv}h_{xy} = 0$. Then
$$
u + \frac{\partial h_{xx}}{\partial u} + 
\frac{\partial h_{xy}}{\partial v}=0. 
$$
Hence we obtain
$$
\frac{\partial (h_{xx}+ \frac{1}{3} u^2) }{\partial u}
=  \frac{\partial (-h_{xy}-\frac{1}{3} uv)}{\partial v}.
\leqno{(2.4)}
$$
Similarly, from the second equality of (2.3), we obtain
$$
\frac{\partial (h_{yy}+ \frac{1}{3} v^2) }{\partial v}
=  \frac{\partial (-h_{xy}-\frac{1}{3} uv)}{\partial u}.
\leqno{(2.5)}
$$
Let $0 < \varepsilon \ll 1$.
For each $p \in \Bbb R^2 = \{(u,v)\}$, let $U_{\varepsilon}(p)$ 
denote the $\varepsilon$-neighborhood of $p$ in 
$\Bbb R^2$.
We now put
$$
(\square_j)_{\varepsilon}: = \bigcup_{p\in \square_j} U_{\varepsilon}(p).
$$ 
Note that (2.4) and (2.5) hold on $\square_j$. To see whether (2.4) and (2.5) 
are true also for 
$(\square_j)_{\varepsilon}$, take an arbitrary point $q$ in $X_j\setminus T$.
If necessary, replace the complex coordinates $(z_1, z_2)$ 
for $T = (\Bbb C^*)^2$ by 
$$(z_1^{\alpha_1}z_2^{\beta_1},z_1^{\alpha_2}z_2^{\beta_2})
\;\;\text{ for some }
\begin{pmatrix}
\alpha_1 &\alpha_2\\
\beta_1 &\beta_2
\end{pmatrix}\; \in\; \operatorname{GL}(2,\Bbb Z).
$$
Then we may assume that $z_1$, $z_2$ regarded as meromorphic functions on 
$X_j$ are holomorphic at $q$ satisfying 
$$
z_1(q) =0.
$$
Put $b:= |z_2(q)|^2 \geq 0$.
Note that there is a real-analytic function $Q = Q(r_1,r_2)$ in two varibles $r_1$, $r_2$
defined in a neighborhood of $(0,b)$ such that $e^{-h}\, =\, Q(|z_1|^2, |z_2|^2)$. Then
$$
h_x = - \,r_1\, Q^{-1} \frac{\partial Q}{\partial r_1}
\quad\text{and}\quad
h_y = - \,r_1\, Q^{-1} \frac{\partial Q}{\partial r_2},
$$
where both are evaluated at $(r_1,r_2) = (|z_1|^2,|z_2|^2)$.
Then by allowing $r_1$ (and $r_2$ as well if $b =0$)
to take negative values (see also [9;p.722-723]), 
we see that the terms $h_{xx}$, $h_{xy}$, $h_{yy}$ are
 well-defined also on $(\square_j)_{\varepsilon}\setminus \square_j$ 
Then both (2.4) and (2.5) hold on $(\square_j)_{\varepsilon}$.
Since $(\square_j)_{\varepsilon}$ is simply connected, 
we now obtain real-analytic functions $H_1 = H_1 (u, v)$ 
and $H_2 = H_2 (u, v)$ on 
$(\square_j)_{\varepsilon}$ such that
$$
\frac{\partial H_1}{\partial v} = h_{xx} + \frac{1}{3}u^2
\quad\text{and}\quad
\frac{\partial H_1}{\partial u} = -h_{xy}-\frac{1}{3}uv;
\leqno{(2.6)}
$$
$$
\frac{\partial H_2}{\partial u} = h_{yy} + \frac{1}{3}v^2
\quad\text{and}\quad
\frac{\partial H_2}{\partial v} = -h_{xy}-\frac{1}{3}uv. 
\leqno{(2.7)}
$$
In view of the second equalities of (2.6) and (2.7), we obtain 
$\frac{\partial H_1}{\partial u} = \frac{\partial H_2}{\partial v}$.
Hence, for some real-analytic function $H=H(u,v)$ 
on $(\square_j)_{\varepsilon}$,
$$
H_1 = \frac{\partial H}{\partial v}\quad\text{and}\quad 
H_2 = \frac{\partial H}{\partial u}.
\leqno{(2.8)}
$$
By (2.8) together with the first equalities of (2.6) and (2.7),
 we obtain the following on $(\square_j)_{\varepsilon}$:
$$
h_{xx} = H_{vv}-\frac{1}{3}u^2, \;\;
h_{yy} = H_{uu}-\frac{1}{3}v^2,\;\;
h_{xy} = -H_{uv}-\frac{1}{3}uv,
\leqno{(2.9)}
$$
where $H_{uu}:= (\partial^2 H)/(\partial u^2)$, 
$H_{uv}:=(\partial^2 H)/(\partial u\partial v)$, 
$H_{vv}:= (\partial^2 H)/(\partial v^2)$.
On the other hand, by (2.2) and (2.9),
\begin{align*}
&\frac{\partial e^{-h}}{\partial u}
= - e^{-h} (h_x \frac{\partial x}{\partial u} + h_y \frac{\partial y}{\partial u})
= - e^{-h} (u h^*_{uu} + v h^*_{uv})\\
&= -uh_{yy} + vh_{xy}= -u(H_{uu}-\frac{1}{3}v^2)+v(-H_{uv}-\frac{1}{3}uv)\\
&=-(uH_{uu}+vH_{uv}) = \frac{\partial}{\partial u}\, (H-uH_u -vH_v).
\end{align*}
Similarly, 
$$
\frac{\partial e^{-h}}{\partial v}=\frac{\partial}{\partial v}\, (H-uH_u -vH_v).
$$\
Hence $\,e^{-h} = H-uH_u -vH_v + C\,$ for some real constant $C$. Replacing $H$ 
by $H+C$, we may assume without loss of generality that
$$
e^{-h} = H-uH_u -vH_v.
\leqno{(2.10)}
$$
In view of (2.9) and (2.10), the equation (2.1) is rewritten as
$$
\vmatrix 
H_{vv}-\frac{1}{3}u^2,& H_{uv}+\frac{1}{3}uv\\
H_{uv}+\frac{1}{3}uv,& H_{uu}-\frac{1}{3}v^2
\endvmatrix
=H-uH_u -vH_v.
\leqno{(2.11)}
$$
Put $s:= H_u$ and $t:= H_v$. In the next section, we shall show that
the image of $\bar{\square}_j$
under the mapping 
$$
\bar{\square}_j \owns (u,v) \mapsto (s(u,v), t(u,v)) \in \Bbb R^2
$$ 
is nothing but $\bar{\Omega}_j$ in the introduction.
Moreover, by this map, the boundary $\partial\square_j$ is mapped
onto the boundary $\partial\Omega_j$.
We now consider the Legendre transform
$\psi := uH_u+vH_v-H$.  Regard $\psi$ as a function in $(s,t)\in \bar{\Omega}_j$.
Since $\operatorname{Hess}(h)$ is positive on $T$ and vanishes on $X_j\setminus T$,
we see from (2.11) that $\psi$ is negative on $\Omega_j$, and
vanishes just on the boundary $\partial\Omega_j$. Then 
$$
\psi_s \,= \,u \quad\text{ and }\quad \psi_t\, =\, v.
\leqno{(2.12)}
$$
Moreover, in view of the equalities $\psi_{ss} = \partial u/\partial s$, 
$\psi_{tt} = \partial v/\partial t$ and $\psi_{st}= \partial u/\partial t$,
we see that
$$
\begin{pmatrix}
\psi_{ss} & \psi_{st}\\
\psi_{st} & \psi_{tt}
\end{pmatrix} \; = \;
 \begin{pmatrix}
H_{uu} & H_{uv}\\
H_{uv} & H_{vv}
\end{pmatrix}^{-1}.
\leqno{(2.13)}
$$
In particular, $\det\operatorname{Hess}(H):= H_{uu}H_{vv}-(H_{uv})^2$ 
and $\det\operatorname{Hess}(\psi):= \psi_{ss}\psi_{tt}-(\psi_{st})^2$ satisfy
$\; \det\operatorname{Hess}(H)\cdot \det\operatorname{Hess}(\psi) =1$.
Now by (2.11),
$$
\det\operatorname{Hess}(H) \; = \; -\psi + \frac{1}{3}(u^2 H_{uu}+2uv H_{uv} + v^2 H_{vv}).
\leqno{(2.14)}
$$
By (2.13), we have $H_{uu} = \psi_{tt}/\det\operatorname{Hess}(\psi)$,
$H_{vv} = \psi_{ss}/\det\operatorname{Hess}(\psi)$, $H_{uv} = -\psi_{st}/\det\operatorname{Hess}(\psi)$. Hence, from (2.12) and (2.14), it follows that
\begin{align*}
1\; &=\;\det\operatorname{Hess}(\psi)
\left\{-\psi + \frac{1}{3}(u^2 H_{uu}+2uv H_{uv} + v^2 H_{vv})\right\}\\
&=\; -\psi\det\operatorname{Hess}(\psi) + \frac{1}{3}
\left\{\psi_s^2\psi_{tt}-2\psi_s\psi_t\psi_{st}+\psi_t^2\psi_{ss}\right\}.
\end{align*}
Thus, we obtain the equality (1.3). By setting (1.2), we finally see that (1.1) holds, 
as required.

\section{The boudary condition}

We first consider the case $j=1$, so that $X_j = \Bbb P^2(\Bbb C)$.
Then the K\"ahler-Einstein form $\omega$ on $X_j$ given by 
$$
h = - \log 9  - x - y + 3\log (1+ e^x +e^y).
$$ 
is known as the Fubini-Study form.
This obviously satisfies the equation (2.1). Moreover, $u := h_x$ and $v :=h_y$ 
satisfy the inequalities
\begin{align*}
&1-(u+v) = \frac{3}{1+e^x+e^y} \geq 0,
\\
&u+1=\frac{3e^x}{1+e^x+e^y}\geq 0,  
\qquad v+1=\frac{3e^y}{1+e^x+e^y}\geq 0.
\end{align*}
In this case, $H$ and $\psi$ are
\begin{align*}
&H = \frac{uv(u+v)}{6} + \frac{u^2+uv+v^2}{3}+\frac{1}{3},\\
&\psi = uH_u + vH_v -H = \frac{1}{3}(u+1)(v+1)(u+v-1) \; \leq \; 0.
\end{align*}
Then $h$ and $H$ satisfy (2.9) and (2.11). Moreover $\psi$, when regarded as a function on $\bar{\square}_j$, is negative
on $\square_j$ vanishing on the boundary $\partial\square_j$.
In addition to this, for $s:= H_u$ and $t:= H_v$, we can easily check that
$$
\begin{cases}
&\rho_{11}(s,t) = 0 \quad\quad\quad\quad \text{on the line $u+v=1$,}\\
&\rho_{21}(s,t) =0 \quad\quad\quad\quad \text{on the line \;\;$u=-1$,}\\
&\rho_{31}(s,t) = 0\quad\quad\quad\quad \text{on the line \;\;$v=-1$,}
\end{cases}
$$
and that the mapping $\bar{\square}_j  \owns (u,v) \mapsto (s(u,v),t(u,v))\in   
\bar{\Omega}_j$ takes $\square_j$ diffeomorphically onto $\Omega_j$. 
We now regard $\psi$ as a function on $\bar{\Omega}_j$. Then $\psi$ is negative
on $\Omega_j$ vanishing on the boundary $\partial\Omega_j$. We also see that $\psi$ is a root of a polynomial of degree 4 with
coefficients in $\Bbb Q [s, t]$ such that the leading coefficient is 1.
Moreover, the asymptotic expansion of $\psi$ along the boundary curve 
$\partial \Omega_j$, especially along $\{u = -1\}\cap\partial\square_j$, shows
$$
\psi\; = \; -\rho_{21}\; +\; \frac{(\rho_{21})^2}{-4s+2t-1}\; + \;
\text{higher order terms in $\rho_{21}$},
$$
where from this expression of $P(\Omega_j )$,
we easily see that the pullback
of the K\"ahler-Einstein form $\omega$ to each (irreducible) component of $X_j \setminus T$ 
is nothing but the Fubini-Study form on $\Bbb P^1 (\Bbb C )$.

\medskip
We next consider the case $j = 2$, so that $X_j = \Bbb P^1 (\Bbb C )\times 
\Bbb P^1 (\Bbb C )$. 
Then we fix the K\"ahler-Einstein form $\omega$ on $X_j$ defined by 
$$
h = - 2\log 2  -  x - y + 2\log (1+ e^x)+ 2 \log (1+e^y).
$$ 
This again satisfies the equation (2.1). Then $u := h_x$ and $v:=h_y$ 
satisfies
\begin{align*} 
&-1\leq \; -1 +\frac{2}{1+e^x} = \; u \;  = 1 - \frac{2 e^x}{1+e^x} \;  \leq 1,\\
&-1 \leq \; -1+\frac{2}{1+e^y} = \; v =\; 1 - \frac{2 e^y}{1+e^y} \;\leq 1.
\end{align*}
In this case, $H$ and $\psi$ are expressible as
\begin{align*}
&H \; =\; \frac{u^2+v^2}{4} - \frac{u^2v^2}{12} + \frac{1}{4},\\
&\psi \; =\; -\frac{1}{4}(1-u^2)(1-v^2)\; \leq \; 0.
\end{align*}
Then $h$ and $H$ satisfy (2.9) and (2.11). Moreover $\psi$, when regarded as 
a function on $\bar{\square}_j$, is negative
on $\square_j$ vanishing on the boundary $\partial\square_j$.
Now for $s:= H_u$ and $s:= H_v$, we see that
$$
\begin{cases}
&\rho_{12}(s,t) = 0\quad\quad \text{ on the line $v=1$,}\\
&\rho_{22}(s,t) =0 \quad\quad \text{ on the line $u=1$,}\\
&\rho_{32}(s,t)=0 \quad\quad \text{ on the line $v=-1$,}\\
&\rho_{42}(s,t) =0\quad\quad \text{ on the line $u=-1$,}
\end{cases}
$$
and the mapping $\bar{\square}_j  \owns (u,v) \mapsto (s(u,v),t(u,v))\in   
\bar{\Omega}_j$ takes $\square_j$ diffeomorphically onto $\Omega_j$. 
Then $\psi$ regarded as a function on $\bar{\Omega}_j$ is negative
on $\Omega_j$ and vanishes on the boundary $\partial\Omega_j$.
We also see that $\psi$ is a root of a polynomial of degree 5 with
coefficients in $\Bbb Q [s, t]$ such that the leading coefficient is 1.
The asymptotic expansion of $\psi$ along the boundary curve is, 
when restricted to a neighborhood of 
$\{u = 1\}\cap\partial\square_j$ 
for instance, 
$$
\psi\; = \; -\rho_{22}\; +\; \frac{(\rho_{22})^2}{6s-2}\; + \;
\text{higher order terms in $\rho_{22}$}.
$$
In view of this expression of $P(\Omega_j )$, it again follows that the pullback
of the K\"ahler-Einstein form $\omega$ to each component of $X_j \setminus T$ 
is the Fubini-Study form on $\Bbb P^1 (\Bbb C )$.

\medskip
We finally consider the case $j =3$, so that 
$X_j = \Bbb P^2(\Bbb C)\# 3 \bar{\Bbb P}^2(\Bbb C )$. 
Then we may assume without loss of generality that the $K$-invariant K\"ahler-Einstein metric $\omega$ on $X_j$
is invariant under the natural action of the subgroup $D_6
\subset \operatorname{Aut}(X_j)$, where $D_6$ denotes
the dihedral group of order 12. 
Hence from the degeneracy condition for the matrix $\operatorname{Hess}(h)$ along 
$X_j \setminus T$, we now see that
$$
H_{|\partial\square_j}\; = \;\frac{1}{6}\,(r +1),
$$
where $r := u^2+uv+v^2$. Put $\tau := (1-u^2)(1-v^2)(1-(u+v)^2)$.
Then the power series expansion of $H$ along the boundary $\partial \square_j$ 
is given by
$$
H \;=\; \frac{1}{6}\,(r +1)\; + \;\frac{\tau}{12(r -3)}\; + \;\sum_{\alpha \geq 2}
\; \eta_{\alpha} \; \tau^{\alpha},
\leqno{(3.1)}
$$
where each coefficient $\eta_{\alpha}= \eta_{\alpha}( r )$ $(\alpha \geq 3)$ is defined inductively from $\eta_2$, $\eta_3$, \dots ,$\eta_{\alpha -1}$ (and their derivatives).
Now by this (3.1), it is easily seen that $s:=H_u(u,v)$ and $t:=H_v(u,v)$ satisfy
$$
\begin{cases}
&\rho_{13}(s,t) = 0\quad\quad \text{ on $\{u+v=1\}\cap \partial\square_j$,}\\
&\rho_{23}(s,t) =0 \quad\quad \text{ on $\{u=1\}\cap \partial\square_j$,}\\
&\rho_{33}(s,t)=0 \quad\quad \text{ on $\{v=-1\}\cap \partial\square_j$,}\\
&\rho_{43}(s,t) =0\quad\quad \text{ on $\{u+v=-1\}\cap \partial\square_j$,}\\
&\rho_{53}(s,t) =0\quad\quad \text{ on $\{u=-1\}\cap \partial\square_j$,}\\
&\rho_{63}(s,t) =0\quad\quad \text{ on $\{v=1\}\cap \partial\square_j$,}\\
\end{cases}
$$
where the map 
$\bar{\square}_j  \owns (u,v) \mapsto (s(u,v),t(u,v))\in \bar{\Omega}_j$
again takes $\square_j$ diffeomorphically onto $\Omega_j$.
Note that the term $\eta_2$ is uniquely determined by 
the ``Fubini-Pick invariant'' $P(\Omega_j )$, and vice versa.
Then by (3.1), we can explicitly describe $h$ (and hence $\omega$) from 
the data $P(\Omega_j )$, since the equalities (2.2), (2.9), (2.10) above allow us 
to recover $h$ from $H$ by
\begin{align*}
&x = \int \frac{(H_{uu}-\frac{1}{3}v^2)du+(H_{uv}+\frac{1}{3}uv)dv}{H-uH_u-vHv},\\
&y = \int \frac{(H_{vv}-\frac{1}{3}u^2)dv+(H_{uv}+\frac{1}{3}uv)du}{H-uH_u-vHv},\\
&h = - \ln (H-uH_u -vH_v).
\end{align*} 
We also see that $\eta_2$ (and hence $P(\Omega_j ))$ is uniquely determined by the pullback of the K\"ahler-Einstein form $\omega$ to $X_j \setminus T$, 
and vice versa. In particular, it is seen that the pullback of the K\"ahler form $\omega$ to 
any irreducible components of $X_j \setminus T$ can never be a K\"ahler-Einstein form. 
The details in this case $j =3$ and also in the case 
$X = \Bbb P^2(\Bbb C )\# 2\bar{\Bbb P}^2 (\Bbb C )$ 
(cf. \S 4) will be published elsewhere.

\section{Concluding remarks} 

The K\"ahler-Ricci soliton on the toric surface
$\Bbb P^2(\Bbb C )\# \bar{\Bbb P}^2 (\Bbb C )$ is explicitly written (cf.\,[7]) 
by solving an ODE. Now the remaining toric surface is 
$X := \Bbb P^2(\Bbb C )\# 2\bar{\Bbb P}^2 (\Bbb C )$.
As in the preceding sections, the equation for a $K$-invariant K\"ahler-Ricci 
soliton on $X$ (cf.\,[16]) is 
$$
\det \operatorname{Hess}(h) = e^{-h-\alpha (u+v)},
\leqno{(4.1)}
$$
where $0\neq \alpha \in \Bbb  R$ is such that $\alpha (u+v)$ is the Hamiltonian function 
for the holomorphic vector field (cf.\,[15]) associated to the K\"ahler-Ricci soliton.
Let $\square := \bar{\square}\setminus \partial \square$ be the interior 
of the polygon $\bar{\square}$ in $\Bbb R^2$ defined by
$$
\bar{\square}\, = \,\{\,(u,v)\in \Bbb R^2\,;\, 
|u| \leq 1, |v|\leq 1, u+v \leq 1\,\}.
$$
Then by the same argument as in obtaining (2.11) from (2.1), we can reduce (4.1)
to the following equation in $H\in C^{\omega}(\bar{\square} )_{\Bbb R}$ 
with a suitable boundary condition: 
\begin{align*}
&\vmatrix 
H_{vv}+2\alpha H_v+\alpha^2 H-\frac{1}{\alpha}u,& H_{uv}+\alpha H_v+\alpha H_u+\alpha^2 H-
\frac{1}{\alpha^2}\\
H_{uv}+\alpha H_v+\alpha H_u+\alpha^2 H-\frac{1}{\alpha^2},& H_{uu}+2\alpha H_u+\alpha^2 H
-\frac{1}{\alpha}v
\endvmatrix  \\
&=H-u(H_u+\alpha H) -v(H_v+\alpha H) + \frac{uv}{\alpha^2}.
\end{align*}
Since $\partial\square$ is a 1-cycle, we have the following compatibility condition
for the boundary values of $H_u$ and $H_v$:
$$
(2-\alpha^2)e^{3\alpha} = 4 e^{2\alpha} - 2(1+\alpha ),
$$
where $\alpha$ is characterized as the nonzero solution of this equation.
Then this fits to the approximate value of the constant $\alpha$ in [6;\,(14)] (see also [15;\,Lemma 2.2]). However, in this K\"ahler-Ricci soliton case, the equation cannot be so simplified as in (1.1) and (1.3).

\medskip
As compared with [5] and [6], the results in this note may give another frame work for numerical studies of K\"ahler-Einstein metrics and K\"ahler-Ricci solitons.
Let me finally remark that parts of this note are in [10], and were announced in Aug., 1987 in the Taniguchi International Symposium at Katata.

\end{document}